\documentclass[11pt, reqno]{amsart}
\usepackage{graphicx}
\usepackage{amssymb}

\newtheorem{theorem}{Theorem}[section]

\theoremstyle{definition}

\def\beq {\begin{equation}}
\def\endq {\end{equation}}

\def\beq {\begin{equation}}
\def\endeq {\end{equation}}

\renewcommand{\epsilon}{\varepsilon}

\begin{document}
\title{Dual Rectangles}
\author{Graham Everest and Jonny Griffiths}
\address{(GE) School of Mathematics, University of East Anglia,
Norwich NR4 7TJ}
\address{(JG) Paston College, Grammar School Road,
North Walsham, NR28 9JL}
\email{g.everest@uea.ac.uk}
\email{jonny.griffiths@ntlworld.com}

\begin{abstract}
This article concerns a notion of duality between rectangles.
A proof is given that only finitely many integral sided pairs of dual rectangles
exist. Then a geometrical group law is shown to hold on the set of all rational
self-dual rectangles. Finally, the arithmetic of a cubic surface is
used to construct new pairs of rational dual rectangles from old, a
technique inspired by the theory of elliptic curves. 
\end{abstract}

\maketitle


\section{Duality}
The concept of {\it duality} in Mathematics is incredibly rich and powerful, occurring in many places.
Even so, it is
not easy to define in general terms. However the essential properties entail that the dual
of the dual is the starting object and that an object can be self-dual.

In this article, we will say that a pair of rectangles, with sides $(a,b)$ and $(c,d)$ are called {\it dual}
if the area of each is the perimeter of the other. We will assume throughout that the rectangles
are lying down: in other words,
that the pairs are written with $a \ge b$ and $c \ge d$.
\begin{theorem}\label{seven}There are precisely seven pairs of dual rectangles with integral sides.
Two are self-dual: $(4,4)$ and $(6,3)$. The remaining five are
$$(6,4) (10,2),\quad (10,3) (13,2),\quad (10,7) (34,1), \quad (13,6) (38,1), \quad (22,5) (54,1).
$$
\end{theorem}

The duality property translates into the following pair of simultaneous
equations:
\begin{equation}\label{basic}
\begin{aligned}ab &= 2c+2d \\ cd &= 2a+2b.
\end{aligned}
\end{equation}
There are infinitely many rational solutions to the equations in (\ref{basic}) as we now show.
Once $b$ and $d$
are fixed there remain two simultaneous equations in the two variables $a$ and $c$.
It turns out these
have a unique solution provided $bd \neq 4$. In this case $a$ and $c$ are given explicitly by
\begin{equation}\label{canda}
a=\frac{2d^2+4b}{bd-4} \mbox{ and } c=\frac{4d+2b^2}{bd-4}.
\end{equation}
Thus any pair of rational values for $b$ and $d$ with $bd\neq 4$ yields
a pair of rational numbers $a$ and $c$ by solving the simultaneous equations.
({\bf Exercise:} what happens when $bd=4$?) What is more, provided $b$ and $d$ are positive and
$bd>4$, the solutions are guaranteed to
consist of positive rational numbers.
It follows that there exist infinitely many pairs of dual rectangles with rational sides.
Theorem~\ref{seven}
asserts that only finitely many of these pairs have integral sides.

\begin{proof}
Eliminating $d$ from the equations in (\ref{basic}) shows that $c$ satisfies a quadratic
equation
\begin{equation}\label{quadratic}2c^2-cab+4(a+b)=0.
\end{equation}
Solving this equation shows that
$$c=\frac{ab\pm \sqrt{a^2b^2-32(a+b)}}{4}.
$$
Therefore, a necessary (but not sufficient) condition for $c$ to be integral is
\begin{equation}\label{quartic}
a^2b^2-32(a+b)=t^2
\end{equation}
for some non-negative integer $t$. Re-arranging (\ref{quartic}) and factorizing yields
\begin{equation}\label{factored}(ab+t)(ab-t)=32(a+b).
\end{equation}
Since $ab+t$ divides $32(a+b)$, it follows that
\begin{equation}\label{inequality}ab+t \le 32(a+b).
\end{equation}
The hypothesis that $b\le a$ implies that the right hand side of (\ref{inequality}) is bounded above
by $64a$. Since $0\le t$, the left hand side of (\ref{inequality}) is bounded below by $ab$. Therefore
$$ab \le 64 a
$$
and cancelling $a$ shows that $b \le 64$. By an entirely symmetrical argument, $d\le 64$
also. Hence there are a finite number of integer pairs $b$ and $d$. Each one yields at most one
integer pair $a$ and $c$ using the same method as before. For each value of $b$ and
$d$ with $1 \le b,d \le 64$ and $bd\neq 4$, the corresponding values
of $a$ and $c$ in (\ref{canda}) can be checked for integrality, using a computer (we did).
The cases when $bd=4$ do not
lead to consistent
equations.
\end{proof}

{\bf Exercise:} Using the formul\ae { }in (\ref{canda}) turns up quite
a few examples of pairs of dual rectangles where three out of the four sides
are integers. For example:
$$(7, 3) (8, 5/2), \quad (7, 5) (16, 3/2), \quad (33, 3) (48, 3/2), \quad (89, 1) (40, 9/2).
$$
Adapt the method of proof above to show there only finitely many pairs of dual rectangles
with three integral sides, and list them all.

\section{Constructing Self-Dual Rectangles}\label{sdrs}
A rational self-dual rectangle corresponds to a solution of
the equation
\begin{equation}\label{selfdual}xy=2x+2y \mbox{ or } (x-2)(y-2)=4,
\end{equation}
in positive rational numbers. The equation (\ref{selfdual}) describes
a rectangular hyperbola.
Suppose we are given two distinct rational points on the hyperbola
$$P=\left(p,\frac{2p}{p-2}\right) \mbox{ and } Q=\left(q,\frac{2q}{q-2}\right) \mbox{ with } p,q >2.
$$
Consider the triangle $OPQ$. The orthocentre of the triangle has rational coefficients.
Explicitly,
\begin{equation}\label{ortho}\left(2+\frac{8}{(p-2)(q-2)},2+\frac{(p-2)(q-2)}{2}\right),
\end{equation}
and this point lies on the hyperbola too (see \cite{scott}) as the following diagram illustrates.

\bigskip
\begin{center}
\includegraphics[height=60mm]{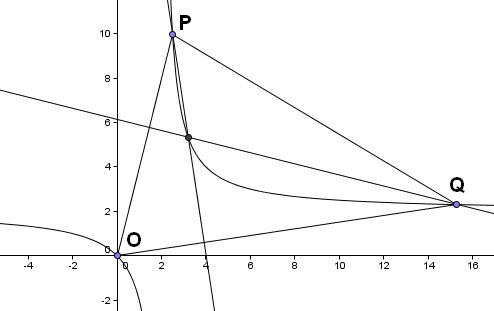}
\end{center}
\bigskip

Now define $P+Q$ as the reflection of the orthocentre of the triangle $OPQ$ in the line $y=x$. Plainly
the point so defined has rational coordinates:
\begin{equation}\label{orthogroupop}\left(p,\frac{2p}{p-2}\right)
+\left(q,\frac{2q}{q-2}\right)=\left(2+\frac{(p-2)(q-2)}{2},2+\frac{8}{(p-2)(q-2)}\right).
\end{equation}
Our remarks confirm that the map defined by (\ref{orthogroupop}) is closed: thus it
is a binary operation. What is more, we can take (\ref{orthogroupop}) to define the
map even when~$p=q$. In other words, we can also `double' points.
\begin{theorem}\label{group}
The operation defined by (\ref{orthogroupop}) is commutative and
associative having the point $(4,4)$ as the identity. Thus the set
of rational self-dual rectangles forms an abelian group.
\end{theorem}

\noindent{\bf Exercise:} prove Theorem \ref{group}.

\section{From Dual Rectangles to Cubic Surfaces}

In this section, we show how the theory of
Diophantine Equations (polynomial equations with integer
coefficients) can be used to construct new solutions to (\ref{basic})
from known solutions. Our motivation comes from a recent
successful attack, by Gerry Tunnell, upon a classical problem about
triangles with rational sides having an integral area (see Koblitz' excellent book~\cite{koblitz}).

The right-angled triangle with sides $3-4-5$ has
area equal to~$6$. The triangle with sides
$$\frac{7}{10} \quad \frac{120}{7} \quad \frac{1201}{70}
$$
is also right-angled and has area equal to~$6$. In a sense
which can be made precise, it is the next simplest such triangle.
In \cite[Section 5.2]{everestward} we used
Tunnell's ideas to point out that there are in
fact infinitely many right-angled triangles
with rational sides having area equal to~$6$. The technique uses the
theory of certain cubic curves called {\it elliptic curves}
(the book \cite{siltate} is an excellent reference
for background reading). The way it turns out, a geometric binary
operation can be defined on the
curve
$$y^2=x^3-36x
$$
which allows known rational solutions to be used to construct other rational solutions.
A dialogue between this curve and triangles then allows other triangles to
be constructed.

One of the fascinating aspects of elliptic curves is the way that
combining quite simple solutions can yield very complicated solutions. We will
model this idea now using the equation in (\ref{quadratic}) and show that a
similar phenomenon is at work.
This equation defines a cubic surface in $3$-space so straight lines should meet it in $3$ points.

\bigskip
\begin{center}
\includegraphics[height=60mm]{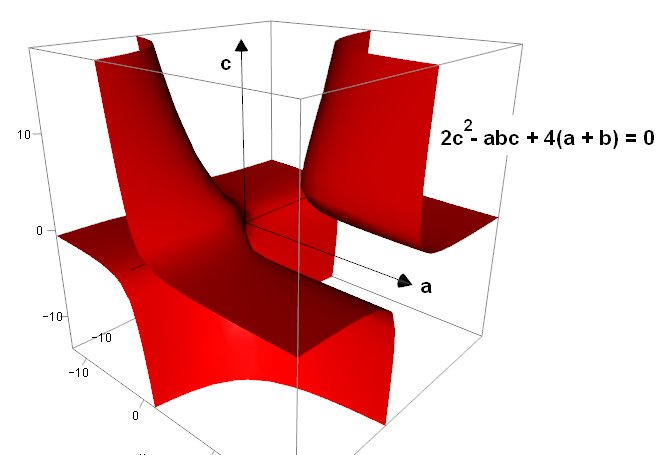}
\end{center}
\bigskip

Define an operation as follows. Given two rational points on
(\ref{quadratic}), join them by a straight line. The third point of intersection
will have rational coordinates. Thus a new solution of (\ref{quadratic}) will be
constructed. Hopefully this will yield a new pair of rational
dual rectangles.

\noindent{\bf Examples}

1. From the list of pairs of integral dual rectangles in Theorem~\ref{seven},
select the rectangle $(6,4)$ which yields the point $(6,4,10)$ on (\ref{quadratic}),
as well as
the rectangle $(22,5)$ which yields the point $(22,5,54)$ on (\ref{quadratic}). The
straight line joining these points has equation
$$(a,b,c)=\theta (6,4,10) + (1-\theta)(22,5,54)
$$
as $\theta$ runs over the real numbers. Substitute into the equation (\ref{quadratic})
and a cubic equation for $\theta$ emerges
\begin{equation}\label{theta}
88\theta^3 -185\theta^2+97\theta=0.
\end{equation}
The $3$ roots of the equation (\ref{theta}) are $0,1$ and $97/88$. The first
two correspond to the points we know about. The last one yields the point
$(48/11,343/88,11/2)$
and this gives rise to the following pair of dual rectangles with
rational sides
$$(48/11, 343/88) (11/2, 727/242).
$$
Notice the way that combining two integral points produced quite a complicated rational point. This
resonates strongly with the theory of elliptic curves. The next example gives
another illustration of the phenomenon.

2. Using the rectangles $(10,3)$ and $(13,6)$ yields the points $(10,3,13)$ and
$(13,6,38)$ on (\ref{quadratic}). The
straight line joining these points has equation
$$(a,b,c)=\theta (10,3,13) + (1-\theta)(13,6,38)
$$
Substituting into the equation (\ref{quadratic})
yields the following cubic equation for $\theta$:
\begin{equation}\label{theta2}
225\theta^3 - 517\theta^2 + 292\theta=0.
\end{equation}
The $3$ roots of the equation (\ref{theta2}) are $0,1$ and $292/225$.
The last one yields the point $(683/75,158/75,50/9)$
and this gives rise to the following pair of dual rectangles with
rational sides
$$(683/75, 158/75) (50/9,2523/625).
$$

{\bf Warning:} This method is not guaranteed to work in a straightforward manner. The geometry of the
cubic surface (\ref{quadratic}) means that the third point of
intersection will not necessarily produce a pair of dual rectangles. If the third point of intersection
exists and has a positive value of~$c$ then the corresponding value of~$d$ exists and
is positive - hence a pair of dual rectangles are created. However a negative or even a zero~$c$ value
could emerge.
In this case the equations make sense but will not allow an interpretation in terms
of dual rectangles.

{\bf Non-Example:} Suppose we try to `add' the two points $(6,4,10)$ and $(10,3,13)$ which
come from Theorem~\ref{seven}.
The
straight line joining these points has equation
$$(a,b,c)=\theta (6,4,10) + (1-\theta)(10,3,13)
$$
as $\theta$ runs over the real numbers. Substitute into the equation (\ref{quadratic})
and a cubic equation for $\theta$ emerges with roots $0,1,13/3$. The third root yields
the point $(-22/3,22/3,0)$ and hence a value for $d=-242/9$. These values satisfy the equations
in (\ref{basic}) but they plainly do not produce a pair of dual rectangles.

\end{document}